\documentclass[11pt]
{amsart}
\usepackage{hyperref}
\hypersetup{nesting=true,debug=true,naturalnames=true}
\usepackage{graphicx,amssymb,upref}
\usepackage{graphicx}
\title{Random walk and Fibonacci matrices}
\author{Theo van Uem}
\address{Theo van Uem\\
Amsterdam University of Applied Sciences, Amsterdam, The Netherlands}
\email{tjvanuem@gmail.com}
\usepackage{amssymb,url}
\newcommand{\NN}{\mathbb N}
\newcommand{\ZZ}{\mathbb Z}
\theoremstyle{plain}
\newtheorem{theorem}[subsection]{Theorem}
\newtheorem{lemma}[subsection]{Lemma}

\newtheorem{proposition}[subsection]{Proposition}
\theoremstyle{definition}

\newtheorem{definition}
{Definition}
\newtheorem{notation}{Notation}
\newtheorem{theo}[subsubsection]{Theorem}
\newtheorem{cor}[subsubsection]{Corollary}
\begin{document}
\begin{abstract}
We study a discrete random walk on a one-dimensional finite lattice, where each state has different probabilities to move one step forward, backward, staying for a moment or being absorbed. We obtain expected number of arrivals and expected time until absorption using a new concept: Fibonacci matrices.
\end{abstract}
\maketitle 
\section{Introduction}
A discrete random walk with variable absorbing probabilities  is described in every state $i \ (i=0,1, \dots ,N)$ by the one-step forward probability $p_i$ , the one-step backward probability $q_i$ , the probability to stay for a moment in the same position $r_i$ and $s_i$ is the probability of absorption in state $i$ where $p_i+q_i+r_i+s_i=1 \ (i=0,1,\dots ,N).$
For this type of random walk we use the notation $[pqrs]$.
In literature (see references) there is a focus on random walks with one or two reflecting and/or absorbing barriers.
 In this paper we have the  freedom of \textit{absorption/reflection in any point at any time with state dependent probabilities}. In this way we can model more complicated situations in physics and operations research. In section 2 we analyze a set of difference equations which is strongly related to the expected number of arrivals and expected time until absorption. Fibonacci numbers and Fibonacci matrices ( a new concept) play an important role in this setting.
In sections 3 and 4 we obtain results for expected number of arrivals and expected time until absorption for a $[pqrs]$  random walk on $[0,N]$. In section 5 we analyze two simple  random walks and their relations to Fibonacci numbers.

\section{Difference equations and Fibonacci matrices}
In section 3 we calculate the expected number of arrivals and we have to solve equations  (see \eqref{5}):
\begin{equation*}
x_n=p_{n-1}x_{n-1}+q_{n+1}x_{n+1}+r_nx_n+\delta (n,i_0)\quad (0\leq n\leq N)
\end{equation*}
and in section 4 we obtain results for expected time until absorption, where we have to deal with (see \eqref{18}):
\[m_i=p_im_{i+1}+q_{i}m_{i-1}+r_{i}m_{i}+1  \ \ (1\leq i\leq N-1)
\]
We shall see that both sets of equations can be handled by solving the next set, which will be the object of research in this section:
\begin{equation}\label{1}
x_{i+1}=\lambda_{m+i}x_i+\mu_{m+i-1}x_{i-1} \ (i=1,2,\dots,N);\ x_0=1;\  x_1=\lambda_m \ (m\in \ZZ)
\end{equation}
We will define Fibonacci matrices which generates in a natural way a unique  solution of the difference equations \eqref{1}.
We start with a Fibonacci sequence:
\begin{equation*}
f_0=f_1=1,\ f_{n+1}=f_n+f_{n-1} \ (n=1,2,\dots)
\end{equation*}
\begin{definition}
 Fibonacci matrices:
$F_0=[1],\ F_1=[\lambda_m] \ (m\in \ZZ),$ where $F_{n+1} \ (n=1,2,\dots)$ with elements $\tau_{ij}^{(m)} \ (i=1,2, \dots,n+1; j=1,2,\dots,f_{n+1})$ is recursively defined by:
\begin{table}[h!]
\centering
%\begin{center}
    \begin{tabular}{|c|c|c|}
    \hline
     $ \tau_{ij}^{(m)}$ & $1$ \dots \dots \dots \dots \dots $f_n$ & $f_n+1$ \dots \dots \dots $f_{n+1}$ \\
    \hline
    $1$&&\\
    .&&\\
     .&$F_n$&$F_{n-1}$\\
      .&&\\
      $n$&&$1$ \dots \dots \dots \dots \dots $1$\\
      $n+1$&$\lambda_{m+n}$ \dots \dots \dots $\lambda_{m+n}$&
      $\mu_{m+n-1}$ \dots \dots $\mu_{m+n-1}$ \\
    \hline
    \end{tabular}%
%\end{center}
\caption{: $F_{n+1}$} 
\label{Table 1}
\end{table}
\end{definition}
so we have:
\begin{equation*}
F_2 = 
\begin{bmatrix}
\lambda_m & 1 \\
\lambda_{m+1} & \mu_m \\
\end{bmatrix}
\end{equation*}

\begin{equation*}
F_3 = 
\begin{bmatrix}
\lambda_m & 1&\lambda_m \\
\lambda_{m+1} & \mu_m &1\\
\lambda_{m+2}&\lambda_{m+2}&\mu_{m+1}\\
\end{bmatrix}
\end{equation*}

\begin{equation*}
F_4 = 
\begin{bmatrix}
\lambda_m & 1&\lambda_m&\lambda_m&1 \\
\lambda_{m+1} & \mu_m &1&\lambda_{m+1} & \mu_m\\
\lambda_{m+2}&\lambda_{m+2}&\mu_{m+1}&1&1\\
\lambda_{m+3}&\lambda_{m+3}&\lambda_{m+3}&\mu_{m+2}&\mu_{m+2}\\
\end{bmatrix}
\end{equation*}
\newline
\begin{lemma}\label{lma1}
$\tau_{i,{f_n+j}}=\tau_{i,j} \quad (1\leq i\leq n-1,\ 1\leq j\leq f_{n-1})$
\end{lemma}
\begin{proof}
The element $F_n$ in Table \ref{Table 1} can be split in $F_{n-1}$ and $F_{n-2}$ (and some $\lambda, \mu$ and $1$ below).  $F_{n-1}$ is in the upper left corner with rows $1$ until $n-1$ and columns $1$ until $f_{n-1}$. The same element   $F_{n-1}$ can be found in  rows $1$ until $n-1$ and columns $f_n+1$ until $f_n+f_{n-1} (=f_{n+1})$.
\end{proof}
\begin{definition}
$F_0^*=1; \  \ \ F_n^*=\sum_{j=1}^{f_n}\prod_{k=1}^n\tau_{kj}^{(m)} \quad (n=1,2,\dots,N+1)$
\end{definition}
\begin{proposition}\label{A}
$F_n^* \ (n=0,1,\dots,N+1)$ is a solution of \eqref{1}
\end{proposition}
\begin{proof}
\begin{equation*}
\begin{aligned}
F_{n+1}^*=
&\sum_{j=1}^{f_{n+1}}\prod_{k=1}^{n+1}\tau_{kj}^{(m)}=
\sum_{j=1}^{f_n}\prod_{k=1}^{n+1}\tau_{kj}^{(m)}+\sum_{j=f_n+1}^{f_{n+1}} \prod_{k=1}^{n+1}\tau_{kj}^{(m)}=\\
&\sum_{j=1}^{f_n}[\prod_{k=1}^{n}\tau_{kj}^{(m)}\tau_{n+1,j}^{(m)}]+ \sum_{j=f_n+1}^{f_{n+1}}[\prod_{k=1}^{n-1}\tau_{kj}^{(m)}\tau_{n,j}^{(m)}\tau_{n+1,j}^{(m)}]=\\
&\lambda_{m+n} \sum_{j=1}^{f_n} \prod_{k=1}^n\tau_{kj}^{(m)}+1.\mu_{m+n-1}\sum_{j=f_n+1}^{f_{n+1}}\prod_{k=1}^{n-1}\tau_{kj}^{(m)}= \\
&\lambda_{m+n} F_n^*+\mu_{m+n-1}\sum_{j=1}^{f_{n-1}} \prod_{k=1}^{n-1}\tau_{kj}^{(m)}=
\lambda_{m+n}F_n^{*}+\mu_{m+n-1}F_{n-1}^{*}
\end{aligned}
\end{equation*}
where we used Lemma ~\ref{lma1} in the penultimate step.
\end{proof} 
\begin{theorem}\label{thm1}
The solution of the linear system \eqref{1} is:
\begin {equation} \label{2}
 x_0=1, \quad x_i=\sum\limits_{j=1}^{f_n}\prod\limits_{k=1}^n\tau_{kj}^{(m)} \quad (i=1,2,\dots,N+1)
\end {equation} 
  where: \\
{\bf Case} $j\leq f_{i+1}$:
\[
\tau{}_{ij}^{(m)}=\left\{\begin{array}{l}\lambda_{m+i-1}\ \ \ \ \ \ \ \    (j=1,2,\dots,f_{i-1}) \\
\mu_{m+i-2}\ \ \ \ \ \  \ \ (j=f_{i-1}+1,\dots,f_i) \\
1\ \ \ \ \ \ \ \ \ \ \ \ \ \ \ \
(j=f_i+1,\dots,f_{i+1})\end{array}\right.
\]
{\bf Case} $j>f_{i+1}$:\\
$\exists n \in \NN$: $1+f_n \leq j \leq f_{n+1}$\\
Let $j_\ell=j-(f_n+f_{n-2}+ \dots +f_{n-2\ell})\quad (\ell=0,1,\dots,
\left \lfloor{\frac{n-1}{2}}\right \rfloor )$ 
and\\
 $k=min\{\ell \in \NN \ | \ j_\ell \leq f_{i+1}\}$,  then: 
\[
\tau{}_{i,j}^{(m)}=\tau_{i,j_k}^{(m)}=\left\{\begin{array}{l}\lambda_{m+i-1}\ \ \ \ \ \ \ \    (j_k=1,2,\dots,f_{i-1}) \\
\mu_{m+i-2}\ \ \ \ \ \  \ \ (j_k=f_{i-1}+1,\dots,f_i) \\
1\ \ \ \ \ \ \ \ \ \ \ \ \ \ \ \
(j_k=f_i+1,\dots,f_{i+1})\end{array}\right.
\]
\end{theorem}

\begin{proof} 
We start with:\\
{\bf Case} $j\leq f_{i+1}$: \\
Substituting \eqref{2} in \eqref{1} yields: \\
\[
\sum\limits_{j=1}^{f_{i+1}}\prod\limits_{k=1}^{i+1}\tau_{kj}^{(m)}=
\lambda_{m+i} \sum\limits_{j=1}^{f_i} \prod\limits_{k=1}^i\tau_{kj}^{(m)}+\mu_{m+i-1}\sum\limits_{j=1}^{f_{i-1}}\prod\limits_{k=1}^{i-1}\tau_{kj}^{(m)}
\]  so:
\[
\sum\limits_{j=1}^{f_{i+1}} \tau_{1j}^{(m)}\tau_{2j}^{(m)}\dots\tau_{i+1,j}^{(m)}= 
\sum\limits_{j=1}^{f_{i}} \tau_{1j}^{(m)}\tau_{2j}^{(m)}\dots\tau_{i,j}^{(m)} \lambda_{m+i} +
\]
\[
\sum\limits_{j=1}^{f_{i-1}} \tau_{1j}^{(m)}\tau_{2j}^{(m)}\dots\tau_{i-1,j}^{(m)}.1.\mu_{m+i-1}
\]
  Using Lemma \eqref{lma1} we get:
  \[
\sum\limits_{j=1}^{f_{i+1}} \tau_{1j}^{(m)}\tau_{2j}^{(m)}\dots\tau_{i+1,j}^{(m)}= 
\sum\limits_{j=1}^{f_{i}} \tau_{1j}^{(m)}\tau_{2j}^{(m)}\dots\tau_{i,j}^{(m)} \lambda_{m+i}+ \]
\[
\sum\limits_{j=f_i+1}^{f_{i+1}} \tau_{1j}^{(m)}\tau_{2j}^{(m)}\dots\tau_{i-1,j}^{(m)}.1.\mu_{m+i-1}.
  \]
  It follows:
  \[\tau_{i+1,j}^{(m)}= \lambda_{m+i} \quad (j\leq f_i)\] and
   \[ \tau_{ij}^{(m)}=1 , \tau_{i+1,j}^{(m)}=\mu_{m+i-1} \quad (f_i+1\leq j\leq f_{i+1}).\]
   {\bf Case} $j>f_{i+1}:$
   \[\exists n \in \NN: 1+f_n \leq j \leq f_{n+1}\] so:
   \[1\leq j_0=j-f_n\leq f_{n-1}\]
   Let $ j_\ell \ \ (\ell=1,2,\dots,
\left \lfloor{\frac{n-1}{2}}\right \rfloor )$ be defined by:
  \[j_\ell=j_{\ell-1}-f_{n-2\ell}\leq f_{n-2\ell-1} \quad \]   
It follows:  \[j_\ell=j-(f_n+f_{n-2}+ \dots +f_{n-2\ell})\quad (\ell=0,1,\dots,
\left \lfloor{\frac{n-1}{2}}\right \rfloor )\] 
Using Lemma \ref{lma1}: 
   $ \tau_{i,j}^{(m)}=\tau_{i,j_0}^{(m)}= \tau_{i,j_1}^{(m)}=\dots=\tau_{i,j_k}^{(m)}$ where $i\leq n-1$.\\
   Let $k=$min$\{\ell \in \NN \ | \ j_\ell \leq f_{i+1}\}$\ then $j_k\leq f_{i+1}$, so we can apply the first part of this proof.
\end{proof}

\begin{notation}
\[
A_i^{(m)}=A_i^{(m)}(\lambda, \mu)=\sum\limits_{j=1}^{f_{i}}\prod\limits_{k=1}^{i}\tau_{kj}^{(m)} \quad  (i=1,2,\dots,N+1).
\]
\end{notation}
where $\lambda=(\lambda_{m},\lambda_{m+1},\dots,\lambda_{m+N}), \quad \mu=(\mu_{m},\mu_{m+1},\dots,\mu_{m+N-1})$ \begin{definition}
 $A_0^{(m)}=1, \quad A_{-1}^{(m)}=0 \quad (m \in \ZZ).$ 
\end{definition}
\begin{theorem}\label{thm2}
$
A_{N+1}^{(m)}=\lambda_{m}A_N^{(m+1)}+\mu_{m} A_{N-1}^{(m+2)}
$
$\quad (m \in \ZZ, \ N=0,1,2,\dots)$
\end{theorem}
\begin{proof}
We write the linear system 
\[x_{i+1}=\lambda_{m+i}x_i+\mu_{m+i-1}x_{i-1} \ (i=1,2,\dots,N);\ x_0=1;\  x_1=\lambda_m \] in matrix notation:
\[
\begin{bmatrix}
1 & 0 & 0 &\dots& \dots&\dots & 0 \\
\lambda_{m} & -1 & 0 & \dots &\dots&\dots & 0 \\
\mu_{m} & \lambda_{m+1} & -1 & \dots&\dots&\dots & 0 \\
\dots  & \dots  & \dots  & \dots &\dots& \dots&\dots  \\
0 & 0&0&\dots & \mu_{m+N-1} & \lambda_{m+N} & -1 \\
\end{bmatrix}
\begin{bmatrix}
x_0 \\ x_1 \\x_2\\ \dots \\ x_{N+1}
\end{bmatrix}
=
\begin{bmatrix}
1 \\0 \\0\\ \dots \\ 0
\end{bmatrix}
\]

The determinant of the matrix is $(-1)^{N+1}$.
Using Cramer's rule we get:
\[
A_{N+1}^{(m)}=x_{N+1}=
\]
\begin{equation*}
(-1)^{N+1}det
\begin{bmatrix}
1 & 0 & 0 &\dots& \dots&\dots &1 \\
\lambda_{m} & -1 & 0 & \dots &\dots&\dots & 0 \\
\mu_{m} & \lambda_{m+1} & -1 & \dots&\dots&\dots & 0 \\
\dots  & \dots  & \dots  & \dots &\dots& \dots&\dots  \\
0 & 0&0&\dots & \mu_{m+N-1} & \lambda_{m+N} & 0 \\
\end{bmatrix}=
\end{equation*}
\begin{equation*}
(-1)^{N+1}(-1)^{N+1}det
\begin{bmatrix}
\lambda_{m} & -1  \\
\mu_{m} & \lambda_{m+1} & -1  \\
  \\
  \\
 &&&&&& \mu_{m+N-1} & \lambda_{m+N}  \\
\end{bmatrix}=
\end{equation*}
\begin{equation*}
\lambda_{m}det
\begin{bmatrix}
\lambda_{m+1} & -1  \\
  \\
  \\
 &&&& \mu_{m+N-1} & \lambda_{m+N}  \\
\end{bmatrix}+
\end{equation*}
\begin{equation*}
\mu_{m}det
\begin{bmatrix}
\lambda_{m+2} & -1  \\
  \\
  \\
 &&&& \mu_{m+N-1} & \lambda_{m+N}  \\
\end{bmatrix}=
\end{equation*}
\[
\lambda_{m}A_N^{(m+1)}+\mu_{m} A_{N-1}^{(m+2)}
\]
\end{proof}
\section{Expected number of arrivals }
We start with two definitions:
\begin{definition}
$p_{ij}^{(k)}$=P(system is in state $j$ after $k$ steps $|$ start in $i$).
\end{definition}
\begin{definition}
$x_j=x_{i,j}=\sum\limits_{k=0}^\infty p_{ij}^{(k)}$
\end{definition}
$x_j$ is the  expected number of arrivals in $j$ when starting in $i$.
We analyze a finite discrete random walk with different absorbing probabilities: in every state $i \ (i=0,1, \dots ,N)$ we have one-step forward probability $p_i$,  one-step backward probability $q_i$, probability to stay for a moment in the same position $r_i$ and $s_i$ is the probability of absorption in state $i$ where $p_i+q_i+r_i+s_i=1 \ (i=0,1,\dots ,N)$ and $q_0=p_N=0$.
The starting point of the random walk on $[0,N]$ is $i_0 \ \  (0\leq i_0\leq N)$.
\begin{theorem}
\begin{equation}\label{5}
x_n=p_{n-1}x_{n-1}+q_{n+1}x_{n+1}+r_nx_n+\delta (n,i_0)\quad (0\leq n\leq N)
\end{equation}
\begin{proof}
When $0< i_0< N$ we have:
\[
x_n=\sum\limits_{k=0}^\infty p_{i_0,n}^{(k)}=p_{i_0,n}^{(0)}+\sum\limits_{k=1}^\infty\sum_{l} p_{i_0,l}^{(k-1)}p_{l,n}=\delta(i_0,n)+\sum_{l} p_{l,n}\sum\limits_{k=1}^\infty p_{i_0,l}^{(k-1)}=
\]
\[\delta(i_0,n)+p_{n-1}x_{n-1}+q_{n+1}x_{n+1}+r_{n}x_{n}
\]
When $i_0=0$ or $i_0=N$, the prove goes along the same lines.
\end{proof}
\end{theorem}

\begin{lemma}\label{98}
The linear system
\begin{equation}\label{99}
x_{i+1}=\lambda_{i} x_i+\mu_{i-1}x_{i-1} \quad (i=i_0+1,i_0+2,\dots,N)
\end{equation}
given $x_{i_0}$ and $x_{i_0+1}$ has solution:
\begin{equation}\label{100}
x_{i_0+k+1}=x_{i_0+1}A_k^{(i_0+1)}(\lambda,\mu)+\mu_{i_0}x_{i_0}A_{k-1}^{(i_0+2)}(\lambda,\mu)\quad (k=1,2,\dots,N-i_0)
\end{equation}
\begin{proof}
We use induction. \\
(i) Substituting $k=1$ in \eqref{100} gives:
\[
x_{i_0+2}=x_{i_0+1}A_1^{(i_0+1)}(\lambda,\mu)+\mu_{i_0}x_{i_0}A_{0}^{(i_0+2)}(\lambda,\mu)=\lambda_{i_0+1}x_{i_0+1}+\mu_{i_0}x_{i_0}
\]
(ii) We give that \eqref{99} is correct for $i=i_0+k-1$ and $i=i_0+k$. We have:
\[
x_{i_0+k+1}=
\lambda_{i_0+k} x_{i_0+k}+\mu_{i_0+k-1}x_{i_0+k-1}=
\]
\[
\lambda_{i_0+k}[x_{i_0+1}A_{k-1}^{(i_0+1)}+\mu_{i_0}x_{i_0}A_{k-2}^{(i_0+2)} ]
+
\mu_{i_0+k-1}[ 
x_{i_0+1}A_{k-2}^{(i_0+1)}+\mu_{i_0}x_{i_0}A_{k-3}^{(i_0+2)}]=
\]
\[
[\lambda_{i_0+k}A_{k-1}^{(i_0+1)}+\mu_{i_0+k-1}A_{k-2}^{(i_0+1)} ]x_{i_0+1}
+
[\lambda_{i_0+k}A_{k-2}^{(i_0+2)}+\mu_{i_0+k-1}A_{k-3}^{(i_0+2)}]\mu_{i_0}x_{i_0}=
\]\
\[
x_{i_0+1}A_k^{(i_0+1)}+\mu_{i_0}x_{i_0}A_{k-1}^{(i_0+2)}
\]
where in the last step we used that $A_i^{(m)}$ is a solution of \eqref{1}; substitute $m=i_0+1,\ i=k-1$ and $m=i_0+2,\ i=k-2$ in \eqref{1}.
\\
So \eqref{99} is also correct for $i=i_0+k+1$.\\
(iii) Apply induction.
\end{proof}
\end{lemma}

In this section we use 
$
A_i^{(m)}(\lambda, \mu)\ (m \geq 0)
$
where 
\[
\lambda=(\lambda_{m},\lambda_{m+1},\dots,\lambda_{N}), \quad \mu=(\mu_{m},\mu_{m+1},\dots,\mu_{N-1})
\]
\[
\lambda_{i}=\frac{1-r_i}{q_{i+1}},\quad 
 \mu_{i}=-\frac{p_{i}}{q_{i+2}}
\]
We also use $
A_i^{(m)}(\rho, \theta) \ (m<0)
$
where 
\[\rho=(\rho_{m},\rho_{m+1},\dots,\rho_{0}), \quad \theta=(\theta_{m},\theta_{m+1},\dots,\theta_{-1})\]
\[
\rho_{-i}=\frac{1-r_i}{p_{i-1}}, \quad 
 \theta_{-i}=-\frac{q_{i}}{p_{i-2}}
\]
\begin{theorem}
{\bf Case} $0<i_0<N$.
\begin{equation}\label{13}
x_{i_0}=[1-r_{i_0}+p_{i_0-1}\theta_{-i_0}\frac{A_{i_0-1}^{(2-i_0)}(\rho,\theta)}{A_{i_0}^{(1-i_0)}(\rho,\theta)}+q_{i_0+1}\mu_{i_0}\frac{A_{N-i_0-1}^{(i_0+2)}(\lambda,\mu)}{A_{N-i_0}^{(i_0+1)}(\lambda,\mu)} ]^{-1}
\end{equation}
For $k=0,1,\dots,N-i_0$ :
\begin{equation}\label{14}
x_{i_0+k+1}=
\end{equation}
\[
\mu_{i_0}x_{i_0}[A_{k-1}^{(i_0+2)}(\lambda,\mu)A_{N-i_0}^{(i_0+1)}(\lambda,\mu)-A_{N-i_0-1}^{(i_0+2)}(\lambda,\mu)A_k^{(i_0+1)}(\lambda,\mu)][A_{N-i_0}^{(i_0+1)}(\lambda,\mu)]^{-1}
\]
and for $k=0,1,\dots,i_0-1$ :
\begin{equation}\label{15}
x_{i_0-(k+1)}=
\end{equation}
\[
\theta_{-i_0}x_{i_0}[A_{k-1}^{(2-i_0)}(\rho,\theta)A_{i_0}^{(1-i_0)}(\rho,\theta)-A_{i_0-1}^{(2-i_0)}(\rho,\theta)A_k^{(1-i_0)}(\rho,\theta)][A_{i_0}^{(1-i_0)}(\rho,\theta)]^{-1}
\]
{\bf Case} $i_0=0$
\begin{equation}\label{17}
x_0=\frac{A_{N}^{(1)}}{{q_1}A_{N+1}^{(0)}}
\end{equation}
\begin{equation}\label{171}
x_i=\frac{A_{N}^{(1)}A_i^{(0)}-A_{i-1}^{(1)}A_{N+1}^{(0)}}{{q_1}A_{N+1}^{(0)}} \quad (i=1,\dots,N)
\end{equation}
\begin{proof}
{\bf Case} $0< i_0< N$.\\
We introduce two artificial states, $N+1$ and $-1$, with specifications:
\begin {equation}\label {3}
x_{N+1}=0; \ \ x_{-1}=0
\end {equation}
\begin {equation}\label {4}
q_{N+1}>0; \ \ p_{-1}>0
\end {equation}
Using \eqref{5} we get
 forward and backward equations:
\begin{equation}\label{6}
x_{i+1}=\frac{1-r_i}{q_{i+1}}x_i-\frac{p_{i-1}}{q_{i+1}}x_{i-1}=\lambda_{i} x_i+\mu_{i-1}x_{i-1} \quad (i=i_0+1,i_0+2,\dots,N)
\end{equation}
\begin{equation}\label{7}
x_{i-1}=\frac{1-r_i}{p_{i-1}}x_i-\frac{q_{i+1}}{p_{i-1}}x_{i+1}=\rho_{-i} x_i+\theta_{-(i+1)}x_{i+1} \quad (i=i_0-1,i_0-2,\dots,0)
\end{equation}
By induction we can prove that solutions of \eqref{6} and \eqref{7} are (see Lemma \eqref{98} for a proof of the first linear system; the second one can be proved the same way):
\begin{equation}\label{8}
x_{i_0+k+1}=x_{i_0+1}A_k^{(i_0+1)}(\lambda,\mu)+\mu_{i_0}x_{i_0}A_{k-1}^{(i_0+2)}(\lambda,\mu)\quad (k=1,2,\dots,N-i_0)
\end{equation}
\begin{equation}\label{9}
x_{i_0-(k+1)}=x_{i_0-1}A_k^{(1-i_0)}(\rho,\theta)+\theta_{-i_0}x_{i_0}A_{k-1}^{(2-i_0)}(\rho,\theta)\quad (k=1,2,\dots,i_0)
\end{equation}
Using \eqref{3},\ \eqref{4},\ \eqref{8} and \eqref{9} we get:
\begin{equation}\label{10}
x_{N+1}=x_{i_0+1}A_{N-i_0}^{(i_0+1)}(\lambda,\mu)+\mu_{i_0}x_{i_0}A_{N-i_0-1}^{(i_0+2)}(\lambda,\mu)=0
\end{equation}
\begin{equation}\label{11}
x_{-1}=x_{i_0-1}A_{i_0}^{(1-i_0)}(\rho,\theta)+\theta_{-i_0}x_{i_0}A_{i_0-1}^{(2-i_0)}(\rho,\theta)=0
\end{equation}
Substituting $n=i_0$ in  \eqref{5} gives :
\begin{equation}\label{12}
(1-r_{i_0})x_{i_0}=1+p_{i_0-1}x_{i_0-1}+q_{i_0+1}x_{i_0+1}
\end{equation}
Using \eqref{10},\ \eqref{11} and \eqref{12}, we get the expected number of arrivals in the starting point $i_0$; see \eqref{13}.\\
For $k=0,1,\dots,N-i_0$ we find the expected number of arrivals in $i_0+k+1$ (use \eqref{12},\ \eqref {10} and \eqref{8}); see \eqref{14}.\\
For $k=0,1,\dots,i_0$  we get \eqref{15} (use \eqref{12},\ \eqref{11} and \eqref{9}).\\
{\bf Case} $i_0=0$ \ ($i_0=N$ proceeds along the same lines).\\
Instead of two artificial states we now need one artificial state $N+1$ with $x_{N+1}=0, \ q_{N+1}>0$.
We get (use \eqref{5}):
\begin{equation*}\label{16}
x_{i+1}=\frac{1-r_i}{q_{i+1}}x_i-\frac{p_{i-1}}{q_{i+1}}x_{i-1}=\lambda_{i} x_i+\mu_{i-1}x_{i-1} \quad (i=1,2,\dots,N)
\end{equation*}
\begin{equation*}
x_{1}=\frac{1-r_0}{q_{1}}x_0-\frac{1}{q_{1}}=\lambda_{0} x_0-\frac{1}{q_{1}} \quad (i=0)
\end{equation*}
with solution (proved by induction):
\[
x_i=x_0A_i^{(0)}-\frac{1}{q_1}A_{i-1}^{(1)} \ (i=1,2,\dots,N+1)
\]
where $A_i^{(m)}=A_i^{(m)}(\lambda,\mu)$.
Using the artificial state we get:
\[
x_{N+1}=x_0A_{N+1}^{(0)}-\frac{1}{q_1}A_{N}^{(1)}=0
\]
resulting in  \eqref{17} and \eqref{171}.
\end{proof}
\end{theorem}
{\bf Remark}.\ \eqref{17} can also be derived from \eqref{13}: use Theorem \eqref{thm2} and $i_0=0, \ p_{-1}=0$.

\section{Expected time until absorption }
Let $T_i$ be the time until absorption when starting in $i \ \ (i=0,1,\dots,N)$.
\begin{definition}
$m_{i}=E[T_i]=\sum_{k=1}^{\infty}kP(T_i=k)\quad (i=0,1,\dots,N)$
\end{definition}
$m_i$ is the expected time until absorption when starting in $i$ .\\
In this section we demand: 
$s_i>0 \quad (i=0,1,..,N)
$.
Let $s=min(s_0,s_1,..,s_N)$. Then P(no absorption after n steps)$\leq (1-s)^n$, so absorption will always occur:
$\sum_{k=1}^{\infty}P(T_i=k)=1$.

\begin{theorem}\label{thm5}
\[
m_0=p_0m_{1}+r_{0}m_{0}+1
\]
\begin{equation}\label{18}
m_i=p_im_{i+1}+q_{i}m_{i-1}+r_{i}m_{i}+1  \ \ (1\leq i\leq N-1)
\end{equation}
\[m_N=q_{N}m_{N-1}+r_{N}m_{N}+1 
\]
\end{theorem}
\begin{proof}
We prove \eqref{18}. The rest is going along the same lines.
\[
m_{i}=E[T_i]=\sum_{k=1}^{\infty}kP(T_i=k)=\sum_{k=1}^{\infty}(k-1)P(T_i=k)+\sum_{k=1}^{\infty}P(T_i=k)=
\]
\[
\sum_{k=2}^{\infty}(k-1)\{p_iP(T_{i+1}=k-1)+q_iP(T_{i-1}=k-1)+r_iP(T_{i}=k-1)\}+1=
\]
\[p_im_{i+1}+q_im_{i-1}+r_im_i+1
\]
\end{proof}
Another way to obtain \eqref{18} is by observing the next step of the random walk: with probability $p_i$ we move to state $i+1$ and then our expectation of time until absorption is $m_{i+1}$. But we did one step, so we have to deal with $1+m_{i+1}$. The last term is about absorption in one step:
\begin{equation}
m_i=p_i(1+m_{i+1})+q_{i}(1+m_{i-1})+r_{i}(1+m_{i})+s_i.1  \ \ (1\leq i\leq N-1)
\end{equation}

In this section we use the next abbreviations:
\[\omega_i=\frac{1-r_i}{p_i} \ (i=0,\dots,N-1); \ \omega_N=1-r_N
\]
\[\phi_i=-\frac{q_{i+1}}{p_{i+1}} \ (i=0,\dots,N-2); \ \phi_{N-1}=-q_N
\]
\[\alpha_i=-\frac{1}{p_i} \ (i=0,\dots,N-1); \ \alpha_N=-1
\]

\begin{theorem}
For $0\leq i\leq N:$
\begin{equation}\label{21}
m_i=\sum_{k=1}^{i}A_{i-k}^{(k)}\alpha_{k-1}+
\end {equation}
\[
-\frac{A_{i}^{(0)}[\omega_{N}\sum_{k=1}^{N}A_{N-k}^{(k)}\alpha_{k-1}+\phi_{N-1}\sum_{k=1}^{N-1}A_{N-1-k}^{(k)}\alpha_{k-1}+\alpha_{N}]}{\omega_{N}A_{N}^{(0)}+\phi_{N-1}A_{N-1}^{(0)}}
\]
\begin{proof}

The $N+1$ forward equations are (using Theorem \eqref{thm5}):
\[
m_{1}=\frac{(1-r_0)}{p_0}m_0-\frac{1}{p_0}=\omega_{0}m_{0}+\alpha_{0}
\]
\[
m_{i+1}=\frac{(1-r_i)}{p_i}m_i-\frac{q_i}{p_i}m_{i-1}-\frac{1}{p_i}=\omega_{i}m_{i}+\phi_{i-1}m_{i-1}+\alpha_{i}\ \ \ (i=1,2,\dots,N-1)
\]
\begin{equation}\label {19}
0=(1-r_N)m_N-q_{N}m_{N-1}-1=\omega_{N}m_{N}+\phi_{N-1}m_{N-1}+\alpha_{N}
\end {equation}

By induction (as in Lemma \eqref{98}) we can prove:

\begin{equation}\label{20}
m_{i}=m_{0}A_{i}^{(0)}+\sum_{k=1}^{i}A_{i-k}^{(k)}\alpha_{k-1} \ \ (i=0,1,\dots,N)
\end{equation}
where 
$A_i^{(m)}=A_i^{(m)}(\omega,\phi)$, with $\omega=(\omega_m,\dots,\omega_N)$ and $\phi=(\phi_m,\dots,\phi_{N-1})$.
Substituting \eqref{20} in \eqref{19} we obtain $m_0$ and again using \eqref{20} we get \eqref{21}. 
\end{proof}
\end{theorem}
\section{Random walk and Fibonacci numbers}
In this section we study two simple  random walks and their relation to Fibonacci numbers.
\subsection{Homogeneous transition probabilities}
We first consider a random walk on $[0,N]$ where we have homogeneous transition probabilities and there is no option to stay for a moment in any state:\\
$p_i=p,\ q_i=q,\ r_i=0,\ s_i=s \ \ (i=1,2,\dots,N-1); \  p+q+s=1;$
$p_0=p,\ s_0=1-p,\ q_N=q,\ s_N=1-q.$
We start in 0.
\begin{theo}
\begin {equation}\label {22}
x_0=x_0^{[N]}=\frac{\sum\limits_{k=0}^{\left \lfloor{\frac{N}{2}}\right \rfloor} \binom {N-k}{k}(-pq)^{k} }{\sum\limits_{k=0}^{\left \lfloor{\frac{N+1}{2}}\right \rfloor} \binom{N+1-k} {k}(-pq)^{k}} \quad (p+q+s=1)
\end{equation}
\begin{proof}
Using \eqref{17} and homogeneity (superscripts in $A_i^{(m)}$ can be omitted and $\lambda_i=\lambda=\frac{1}{q}, \ \mu_i=\mu=-\frac{p}{q}$) we get when $N=1$: \begin{equation} \label{24}
x_0^{[1]}=\frac{A_{1}^{(1)}}{{q_1}A_{2}^{(0)}}=\frac{A_{1}}{{q}A_{2}}=\frac{\lambda}{q(\lambda^2+\mu)}=\frac{1}{1-pq}
\end{equation} 
Using \eqref{17} and Theorem \eqref{thm2} yields:
\[
x_0^{[j]}=\frac{A_{j}^{(1)}}{{q_1}A_{j+1}^{(0)}}=\frac{A_{j}}{{q}A_{j+1}}
\]
so:
\[
\frac{A_{j}}{A_{j+1}}=qx_0^{[j]}
\]
\[
x_0^{[j+1]}=
\frac{A_{j+1}}{{q}A_{j+2}}=\frac{A_{j+1}}{q \{\lambda A_{j+1}+\mu A_{j}\} }=\frac{A_{j+1}}{ A_{j+1}-pA_{j} }=
\frac{1}{1-p\frac{A_j}{A_{j+1}}} =
\frac{1}{1-pqx_0^{[j]}}
\]
We get:
\[
x_0^{[2]}=\frac{1-pq}{(1-pq)-pq.1}=\frac{1-pq}{1-2pq}
\]
\[
x_0^{[3]}=\frac{1-2pq}{(1-2pq)-pq(1-pq)}=\frac{1-2pq}{1-3pq+p^2q^2} \quad etc.
\]
This leads to 
 \eqref{22}, which can be proven by using induction to $N$: \\
(i) N=1 is correct: see \eqref{24}. \\
(ii) Suppose \eqref{22} is correct up to $N+1$.
We first rewrite \eqref{22} to 
$x_0^{[N]}=\frac{\sigma_{N}}{\sigma_{N+1}}$, then:
\[
x_0^{[N+1]}=\frac{1}{1-pqx_0^{[N]}}=\frac{\sigma_{N+1}}{\sigma_{N+1}-pq\sigma_N}
\]
and for all terms except the first and last one in $\sigma_{N+1}-pq\sigma_N$:
\[
\sigma_{N+1}-pq\sigma_N=
\sum \binom{N-k+1}{k}(-pq)^k-pq\sum \binom{N-k}{k}(-pq)^k=\]
\[
\sum \{\binom{N-k+1}{k}+\binom {N-k+1}{k-1}\}(-pq)^k=\sum \binom {N-k+2}{k}(-pq)^k=\sigma_{N+2}
\]
The first term in $\sigma_{N+1}-pq\sigma_N$ is the first term in $\sigma_{N+1}$; the last term in $\sigma_{N+1}-pq\sigma_N$ is the last term in $-pq\sigma_N$ with index $ \lfloor{\frac{N}{2}}\rfloor+1=\lfloor{\frac{N+2}{2}}\rfloor$.\\
So \eqref{22} is also correct for $N+2$.\\
(iii) Apply induction to N.
\end{proof}
\end{theo}
\subsection{Partial absorbing barriers in the endpoints}
Our next random walk is more restricted:
We study  simple random walk with partial absorbing barriers in the endpoints. See the previous random walk but now with $s_i=0 \ (i=1,2,\dots,N-1),\quad p+q=1.$
\begin{theo}
\begin{equation}\label{25}
    \sum_{k=0}^{\lfloor \frac{N}{2} \rfloor} \binom{N-k}{k}(-pq)^k=\left\{
                \begin{array}{ll}
              
                  \frac{q^{N+1}-p^{N+1}}{q-p} \ \ (p \neq q)\\
                  \\
                  \frac{N+1}{2^{N}} \ \ \ (p=q=\frac{1}{2})
                \end{array}
              \right.
\end{equation}
\begin{proof}
The well known probability of absorption in state 0 is:
\[
qx_0=\frac{1-(\frac{p}{q})^{N+1}}{1-(\frac{p}{q})^{N+2}}
\quad (p\neq q)
\] 
So we have: 
\[x_0=\frac{q^{N+1}-p^{N+1}}{q^{N+2}-p^{N+2}} \quad (p\neq q)\]
Using \eqref{22} we guess the first part of \eqref{25}, 
which can be proven by induction to $N$.\\
The second part  of \eqref{25} is obtained by applying de l'Hospitals rule on the first result, or by induction to $N$.
\end{proof}
\end{theo}
\begin{theo} \label{26}
\begin{equation} \label{27}
\sum\limits_{k=0}^{\left \lfloor{\frac{N}{2}}\right \rfloor} \binom{N-k}{k}t^{k}=\frac{1}{2^N}\sum\limits_{n=0}^{\left \lfloor{\frac{N}{2}}\right \rfloor} \binom{N+1}{2n+1}(1+4t)^n 
\end{equation}
\begin{proof}
Let $f_{N,t}=\sum\limits_{k=0}^{\left \lfloor{\frac{N}{2}}\right \rfloor} \binom{N-k}{k}t^{k}$ be the continuation of $\sum\limits_{k=0}^{\left \lfloor{\frac{N}{2}}\right \rfloor} \binom{N-k}{k}(-pq)^{k}$ to $\mathbb R$. \\
Taking $t=-pq$ gives $p^2-p-t=0$ and:
\[
f_{N,t}=\frac{q^{N+1}-p^{N+1}}{q-p}=\frac{(1+\sqrt{1+4t})^{N+1}-(1-\sqrt{1+4t})^{N+1}}{2^{N+1}\sqrt{1+4t}}=
\]
\[
\frac{1}{2^N}\sum\limits_{n=0}^{\left \lfloor{\frac{N}{2}}\right \rfloor} \binom{N+1}{2n+1}(1+4t)^n 
\]
\end{proof}
\end{theo}
\begin{cor}
\[
f_{N}=\frac{1}{2^N}\sum\limits_{n=0}^{\left \lfloor{\frac{N}{2}}\right \rfloor} \binom{N+1}{2n+1}.5^n \quad(N=0,1,2,\dots)
\]
\end{cor}
By repeated differentiating of $f_{N,t} $ we get the 'moments' of  $f_{N}=\sum\limits_{k=0}^{\left \lfloor{\frac{N}{2}}\right \rfloor} \binom{N-k}{k}$
e.g.
\[
(f_{N,t}^{'})_{t=1}=\sum\limits_{k=1}^{\left \lfloor{\frac{N}{2}}\right \rfloor}k \binom{N-k}{k}=
\frac{1}{2^{N-2}}\sum\limits_{n=1}^{\left \lfloor{\frac{N}{2}}\right \rfloor} \binom{N+1}{2n+1}.n.5^{n-1}
\]

\end{document}